# ON THE NEYMAN–PEARSON PROBLEM FOR LAW-INVARIANT RISK MEASURES AND ROBUST UTILITY FUNCTIONALS[1]

By Alexander Schied

*Technische Universität Berlin*

Motivated by optimal investment problems in mathematical finance, we consider a variational problem of Neyman–Pearson type for law-invariant robust utility functionals and convex risk measures. Explicit solutions are found for quantile-based coherent risk measures and related utility functionals. Typically, these solutions exhibit a critical phenomenon: If the capital constraint is below some critical value, then the solution will coincide with a classical solution; above this critical value, the solution is a superposition of a classical solution and a less risky or even risk-free investment. For general risk measures and utility functionals, it is shown that there exists a solution that can be written as a deterministic increasing function of the price density.

**1. Introduction.** Suppose an economic agent wishes to raise the capital $v \geq 0$ today by issuing a contingent claim with a fixed maturity. Suppose furthermore that the (discounted) liability at maturity shall be bounded by some constant $K$. There are many ways of constructing such contingent claims; for instance, the agent could just take out a loan of size $v$, which would lead to the certain liability $-v$ at maturity. Here, our goal is to find a contingent claim such that the *risk* of the terminal liability is minimal among all claims satisfying the issuer's capital constraints.

In a mathematical model, the payoff of a contingent claim is usually described as a random variable $X$ on a probability space $(\Omega, \mathcal{F}, P)$, and we assume that the price of $X$ is given by the expectation $E[\varphi X]$, where the price density $\varphi$ is a $P$-a.s. strictly positive random variable with $E[\varphi] = 1$;

Received February 2003; revised July 2003.
[1]Supported by a Heisenberg fellowship and by the DFG Research Center "Mathematics for key technologies" (FZT 86).

for the purpose of this introduction, we will also assume that $\varphi$ has a continuous distribution. The risk of the liability $-X$ will be measured in terms of a certain risk measure $\rho$. Thus, we are interested in the following problem:

$$
(1) \quad \begin{array}{c} \text{minimize } \rho(-X) \text{ under the constraints} \\ \text{that } 0 \leq X \leq K \text{ and } E[\varphi X] \geq v. \end{array}
$$

Problems of this type arise in various contexts such as in statistical test theory for composite alternatives or in the construction of Arrow–Debreu equilibria in mathematical economics; see, for example, [15] and Chapter 3 of [13]. Our original motivation stems from the problem of finding risk-minimizing hedging strategies in dynamic financial markets. It is well known that such an optimal strategy can be constructed by hedging a solution to a static problem of type (1); see [[3], [9], [10], [18], [19] and [22]].

For the choice $\rho(-X) = E[X]$, the solution to (1) is given by the classical Neyman–Pearson lemma, and for this reason we will call our problem (1) the *Neyman–Pearson problem* for the risk measure $\rho$. The case in which $\rho(-X) = -E[u(-X)]$ for a strictly concave utility function $u$ is also standard. In this article, our main goal is to solve (1) for cases in which the simple expectation $E[\cdot]$ in the two preceding examples is replaced by a supremum (or infimum) of expectations, taken over a nontrivial set $\mathcal{Q}$ of absolutely continuous probability measures. Thus, we are interested in risk measures of the form

$$
(2) \quad \rho(-X) = \sup_{Q \in \mathcal{Q}} E_Q[X]
$$

or, for a utility function $u$,

$$
(3) \quad \rho(-X) = -\inf_{Q \in \mathcal{Q}} E_Q[u(-X)].
$$

The choice of (2) is motivated by the theory of coherent measures of risk as initiated by Artzner, Delbaen, Eber and Heath [1] and further developed by Delbaen [6, 7]. Robust utility functionals of the form (3) arise as a robust Savage representation of preferences on payoff profiles and were suggested by Gilboa and Schmeidler [14]. Both approaches can be brought together by introducing the notion of a convex measure of risk [11, 12], an example being

$$
(4) \quad \rho(-X) = \inf \left\{ m \in \mathbb{R} \,\Big|\, \inf_{Q \in \mathcal{Q}} E_Q[u(m - X)] \geq u(0) \right\}.
$$

We will also obtain results for risk measures of this type. We refer to [13] for surveys on robust Savage representations and risk measures, as well as for standard facts on problems like (1).

The study of general Neyman–Pearson problems for risk measures of the form (2) was initiated by Huber and Strassen [16] and recently continued by Cvitanić and Karatzas [4]. Kirch [17] extended the latter results to robust utility functionals. On the one hand, these articles deal with very general settings, in particular with nonlinear pricing rules of the type $X \mapsto \inf_{P^* \in \mathcal{P}} E^*[X]$, and they yield an interpretation of solutions as classical solutions with respect to "least favorable pairs" $\widehat{Q}, \widehat{P}$ (which in [4] and [17] need not be probability measures). On the other hand, these results rely on essentially nonconstructive methods and typically do not yield explicit solutions. Only a few special cases were solved by Österreicher [21], Rieder [23] and Bednarski [2].

Here, our goal is to obtain explicit solutions to (1) and to point out certain critical phenomena that arise as a consequence of taking suprema (or infima) of expectations in (2) and (3). To this end, we consider a more specific setting with the linear pricing rule $X \mapsto E[\varphi X]$ and make the key assumption that the risk measure $\rho$ is *law-invariant* in the sense that $\rho(-X) = \rho(-Y)$ whenever $X$ and $Y$ have the same law under $P$. While this assumption might be somewhat restrictive from the point of view of theoretical economics, it is satisfied for most risk measures used by practitioners and allows for some interesting mathematical structure. It is satisfied, for instance, if the set $\mathcal{Q}$ in (2) and (3) is of the form

$$(5) \qquad \mathcal{Q}_\lambda = \left\{ Q \ll P \,\Big|\, \frac{dQ}{dP} \le \frac{1}{\lambda} \right\}$$

for some given $\lambda \in (0,1]$ (note that $\mathcal{Q}_1 = \{P\}$). In Section 3, we will solve the Neyman–Pearson problem for

$$(6) \qquad \rho_\lambda(-X) := - \min_{Q \in \mathcal{Q}_\lambda} E_Q[u(-X)],$$

where $u : [0, K] \to \mathbb{R}$ is a strictly concave and continuously differentiable utility function. In particular, we will show that there exists a *critical value* $v_\lambda^* \in (0, K)$ such that the solution $X_v^*$ to the Neyman–Pearson problem for $\rho_\lambda$ coincides with the classical solution $Y_v^0$ for $\rho_1$ as long as $v \le v_\lambda$. For $v > v_\lambda$, however, a *diversification effect* occurs: $X_v^*$ is now a superposition of a risk-free loan of size $\beta \in (0, v)$ and a classical solution $Y_v^\beta$ for $\rho_1$ but with modified upper bound $K - \beta$ and price $v - \beta$. Thus, the solution is of the form

$$(7) \qquad X_v^* = \begin{cases} Y_v^0, & \text{for } v \le v_\lambda, \\ \beta + Y_v^\beta, & \text{for } v > v_\lambda. \end{cases}$$

We will see that, intuitively, this effect is related to an "aversion" of the investor to accept risky bets outside a region of the form $\{\varphi > y\}$, so that capital that cannot be raised by issuing a risky bet on high-price scenarios $\omega \in \{\varphi > y\}$ must instead be obtained via a risk-free loan. We also get a similar result for the translation invariant modification (4) of $\rho_\lambda$.

In the case $u(x) = x$, the problem reduces to the Neyman–Pearson problem for the coherent risk measure

$$\text{(8)} \qquad \text{AVaR}_\lambda(-X) = \max_{Q \in \mathcal{Q}_\lambda} E_Q[X],$$

that will be called the *average value at risk*. It is also known as "conditional value at risk" or "expected shortfall," and coincides, for atomless probability spaces, with the worst conditional expectation

$$\text{WCE}_\lambda(-X) = \sup\{E[X|A] \mid P[A] > \lambda\}.$$

$\text{WCE}_\lambda$ was suggested by Artzner, Delbaen, Eber and Heath [1] as a coherent alternative to the practitioner's value at risk. The Neyman–Pearson problem for $\text{AVaR}_\lambda$ is relatively easy and closely related to results in [2] and [23], as will be explained in Remark 4.7. The solution $X_v^*$ is of the same type as (7), with $Y_v^0$ now denoting the optimal statistical test as provided by the classical Neyman–Pearson lemma. Thus, we have $Y_v^0 = K \cdot \mathbb{I}_{\{\varphi > b\}}$, which can be interpreted as a digital option that pays off in high-price scenarios. Moreover, the critical value $v_\lambda$ can be characterized in terms of the distribution of $\varphi$, and it turns out that $Y_v^\beta = (1 - \beta)Y_{v_\lambda}^0$, thus determining $\beta$ as $(v - v_\lambda)/(1 - v_\lambda)$.

This solution for $\text{AVaR}_\lambda$ will be obtained as a corollary to the more general Theorem 4.1. It solves the Neyman–Pearson problem for the class of *quantile-based coherent measures of risk* that was introduced by Kusuoka [20]. Such a risk measure is of the form

$$\text{(9)} \qquad \rho_k(X) = \int_0^1 k(t) q_{-X}(t) \, dt,$$

where $k : [0, 1) \to [0, \infty)$ is an increasing right-continuous function such that $\int_0^1 k(t) \, dt = 1$, and where $q_X$ denotes a quantile function of the random variable $X \in L^\infty$. $\text{AVaR}_\lambda$ corresponds to the choice $k = \frac{1}{\lambda}\mathbb{I}_{[1-\lambda,1)}$. Moreover, Kusuoka [20] showed that all law-invariant coherent risk measures which admit a representation (2) can be constructed from this class of quantile-based coherent risk measures. The maximal representing set $\mathcal{Q}$ for $\rho_k$ has been described by Dana and Carlier [5].

The Neyman–Pearson problem for $\rho_k$ of (9) admits a solution of the form

$$X^* = \beta \cdot \mathbb{I}_{[a,b)}(\varphi) + K \cdot \mathbb{I}_{[b,\infty)}(\varphi),$$

where the parameters $0 \le \beta < K$ and $0 \le a \le b \le \infty$ can be obtained via a nonlinear variational problem, which involves only two real parameters and which can be solved in a straightforward manner. In contrast to the case of $\text{AVaR}_\lambda$, one may encounter the case $0 < a < b < \infty$, which now corresponds to a diversification into the *two digital options* $(K-\beta)\cdot \mathbb{I}_{[b,\infty)}(\varphi)$ and $\beta \cdot \mathbb{I}_{[a,\infty)}(\varphi)$, the latter being less risky than the former but no longer risk-free.

Our method in obtaining these results is different from the ones used by Huber and Strassen [16], Cvitanić and Karatzas [4], Kirch [17] and others. It is based on the key observation that, for a large class of law-invariant risk measures $\rho$, there exists a deterministic increasing function $f^*:(0,\infty) \to [0,K]$ such that $X^* := f^*(\varphi)$ solves (1). Thus, we are able to reduce the original problem for risk measures such as (6) or (9) to a semiclassical problem of Neyman–Pearson type, but with the additional constraint that the solution must be an increasing function of the price density. If $\rho$ involves the set $\mathcal{Q}_\lambda$ of (5), then this auxiliary problem can be solved directly. In the case of a general quantile-based coherent risk measure, the auxiliary problem is first transformed into a moment problem for subprobability measures, which then can be solved by using general integral representation results.

This paper is organized as follows. In Section 2, we will look at general properties of solutions to the Neyman–Pearson problem (1), assuming only that our risk measure satisfies certain "axioms." We will comment on the existence and (non)-uniqueness of solutions, and we will prove our key result on the existence of a deterministic increasing function that yields a solution when applied to the price density $\varphi$. In Section 3, we will solve the Neyman–Pearson problem for robust utility functionals (6) and their translation invariant modification. In Section 4, we will consider quantile-based coherent risk measures of the form (9). In a first step, we will show that solving a simple moment problem within a small class of two-step functions yields also solutions to our Neyman–Pearson problem. In a second step, we further reduce the moment problem to a two-dimensional variational problem. Section 5 contains the proof of the first reduction theorem in Section 4.

**2. The general structure of solutions.** In this section we discuss the general structure of solutions to the Neyman–Pearson problem (1), where we take $\rho$ as a real-valued functional on $L^\infty := L^\infty(\Omega,\mathcal{F},P)$ that satisfies the following properties for all $X,Y \in L^\infty$:

(10) Monotonicity: If $X \leq Y$, then $\rho(-X) \leq \rho(-Y)$.

(11) Convexity: $\rho(\lambda X + (1-\lambda)Y) \leq \lambda \rho(X) + (1-\lambda)\rho(Y)$ for $0 \leq \lambda \leq 1$.

(12) Law invariance: If $X$ and $Y$ have the same law under $P$, then
$$\rho(X) = \rho(Y).$$

For simplicity, we will also assume that

(13) $\mathbb{R} \ni m \longmapsto \rho(-m)$ is continuous and strictly increasing on $[0,K]$.

Clearly, this property holds if $\rho$ satisfies the additional axiom of

(14) Translation invariance: $\rho(X+m) = \rho(X) - m$ for $m \in \mathbb{R}$ and $X \in L^\infty$,

in which case $\rho$ is a law-invariant *convex measure of risk* [1, 11, 13]. We also suppose that $\rho$ is continuous from above:

(15) $\qquad X_n \searrow X, \quad P\text{-a.s.} \quad \Longrightarrow \quad \rho(X_n) \nearrow \rho(X).$

It is straightforward to check that, given the monotonicity of $\rho$, continuity from above is equivalent to the so-called *Fatou property*:

(16) $$\rho(X) \leq \liminf_{n \uparrow \infty} \rho(X_n)$$
$$\text{for all bounded } (X_n)_{n \in \mathbb{N}} \subset L^\infty \text{ with } X_n \longrightarrow X\, P\text{-a.s.};$$

see, for example, Lemma 4.16 in [13]. Standard arguments such as those in Remark 3.39 of [13] then show:

LEMMA 2.1 (Existence of solutions). *Under conditions* (10), (11) *and* (15), *there exists a solution to the Neyman–Pearson problem* (1).

We will also assume throughout this paper that the underlying probability space $(\Omega, \mathcal{F}, P)$ is atomless. This condition guarantees that $\rho$ is defined on a sufficiently large domain, and it is equivalent to the existence of a random variable with a continuous distribution.

REMARK 2.2 (Nonuniqueness of solutions). The solution to the Neyman–Pearson problem need not be unique. Take, for example, $\rho(-X) = E[X]$ and consider the solution

(17) $$X^* = K \cdot \mathbb{I}_{\{\varphi > c\}} + \gamma \cdot \mathbb{I}_{\{\varphi = c\}}$$

for certain constants $c \geq 0$ and $\gamma \in [0, K]$ as provided by the classical Neyman–Pearson problem. If the distribution of $\varphi$ is not continuous at $c$, one typically has $\gamma \in (0, K)$, and the usual randomization of $X^*$ yields another solution $\widetilde{X}$ which takes only the values 0 and 1. More precisely, $\widetilde{X}$ coincides with $X^*$ on $\{\varphi \neq c\}$; otherwise $\widetilde{X}$ is either 0 or 1, according to an independent Bernoulli experiment with success probability $\gamma$. If one insists on $\sigma(\varphi)$-measurable solutions, then (17) is the only such solution. But uniqueness may also fail in the class of $\sigma(\varphi)$-measurable solutions as will be shown in Remark 4.3. On the other hand, uniqueness in the class of $\sigma(\varphi)$-measurable solutions implies uniqueness in the class of $\mathcal{F}$-measurable solutions, provided that the price density has a continuous distribution; see Proposition 2.7.

We continue with the following general lemma that was suggested by Hans Föllmer and that is of independent interest.

LEMMA 2.3. *Let $\mathcal{G} \subset \mathcal{F}$ be a countably generated $\sigma$-algebra. Then*

$$\rho(X) \geq \rho(E[X|\mathcal{G}]) \qquad \text{for all } X \in L^\infty.$$

*In particular,*

(18) $$\rho(X) \geq \rho(E[X]).$$

PROOF. Lemma 4.45 in [13] states that

(19) $$\rho(X) \geq \rho(X \cdot \mathbb{I}_{A^c} + E[X|A] \cdot \mathbb{I}_A)$$

for any set $A \in \mathcal{F}$ with $P[A] > 0$ (note that the proof of the cited lemma does not use the translation invariance of $\rho$). Let $B_1, B_2, \ldots$ be a sequence of sets in $\mathcal{F}$ such that $\mathcal{G} = \sigma(B_1, B_2, \ldots)$, and denote by $A_1, \ldots, A_m$ the atoms in $\mathcal{G}_n := \sigma(B_1, \ldots, B_n)$. Applying (19) successively with $A := A_1, A_2, \ldots, A_m$ yields

$$\rho(X) \geq \rho\left(\sum_{i=1}^m E[X|A_i] \cdot \mathbb{I}_{A_i}\right) = \rho(E[X|\mathcal{G}_n]).$$

Thus, by the martingale convergence theorem and the Fatou property (16),

$$\rho(E[X|\mathcal{G}]) \leq \liminf_{n\uparrow\infty} \rho(E[X|\mathcal{G}_n]) \leq \rho(X).$$

Finally, (18) follows by taking $\mathcal{G} = \{\varnothing, \Omega\}$. □

The first consequence of the preceding lemma is that the price constraint in problem (1) can be reduced to an equality:

LEMMA 2.4. *Any solution $X^*$ of the Neyman–Pearson problem with capital constraint $v \in [0, K]$ satisfies $E[\varphi X^*] = v$.*

PROOF. The case $v \in \{0, K\}$ is trivial, and so it is enough to consider $v \in (0, K)$. Note that (18) implies that any solution $X^*$ satisfies $\rho(-X^*) > \rho(0)$. Indeed, since $E[\varphi X^*] \geq v > 0$ and $X^* \geq 0$, we must have $E[X^*] > 0$, and (18) and (13) yield $\rho(-X^*) > \rho(0)$. Now suppose by way of contradiction that $E[\varphi X^*] > v$. Then we define $\widetilde{X} := \alpha X^*$, where $\alpha := v/E[\varphi X^*] < 1$. The convexity of $\rho$ implies that

$$\rho(-\widetilde{X}) = \rho(-\alpha X^* - (1-\alpha)0)$$
$$\leq \alpha\rho(-X^*) + (1-\alpha)\rho(0)$$
$$< \rho(-X^*),$$

which, in view of $E[\varphi \widetilde{X}] = v$, contradicts the optimality of $X^*$. □

Another immediate consequence of Lemma 2.3 is the following: If $X^*$ solves the Neyman–Pearson problem (1), then so does $\widetilde X^* := E[X^*\,|\,\varphi]$. In particular, there always exists a $\sigma(\varphi)$-measurable solution. The following key proposition states a crucial property of such solutions. Note that we always use the term "increasing function" synonymously to "nondecreasing function."

PROPOSITION 2.5. *Every $\sigma(\varphi)$-measurable solution $X^*$ can be written as $X^* = f^*(\varphi)$ for some deterministic increasing function $f^*$.*

The proof of this proposition is based on the following version of the classical Hardy–Littlewood inequalities, which we recall here for the convenience of the reader. See, for example, Theorem 2.76 of [13] for a proof.

THEOREM 2.6 (Hardy–Littlewood). *Let $X$ and $Y$ be two nonnegative random variables, and let $q_X$ and $q_Y$ denote quantile functions of $X$ and $Y$. Then,*

$$\int_0^1 q_X(1-t)q_Y(t)\,dt \le E[XY] \le \int_0^1 q_X(t)q_Y(t)\,dt.$$

*If $X = f(Y)$, then the lower (upper) bound is attained if and only if $f$ can be chosen as a decreasing (increasing) function.*

We will also need the following property of quantile functions: If $f$ is an increasing function and $Y$ is a nonnegative random variable, then the quantile $q_{f(Y)}$ of $f(Y)$ satisfies

(20) $$q_{f(Y)}(t) = f(q_Y(t)) \qquad \text{for a.e. } t \in (0,1);$$

see, for example, Lemma 2.77 in [13].

PROOF OF PROPOSITION 2.5. Since the underlying probability space is atomless, there exists a random variable $U$ with a uniform distribution on $(0,1)$ such that $\varphi = q_\varphi(U)$. Now let $X^*$ be any solution to the Neyman–Pearson problem [for further application of this argument in Proposition 2.7, we do not yet assume that $X^*$ is $\sigma(\varphi)$-measurable]. Denote by $F_\varphi$ the distribution function of $\varphi$, and define

$$f(x) = \begin{cases} q_{X^*}(F_\varphi(x)), & \text{if } F_\varphi \text{ is continuous at } x, \\ \dfrac{1}{F_\varphi(x) - F_\varphi(x-)} \displaystyle\int_{F_\varphi(x-)}^{F_\varphi(x)} q_{X^*}(t)\,dt, & \text{otherwise.} \end{cases}$$

Then $f$ is increasing, and $X := f(\varphi)$ satisfies

(21) $$X = E[q_{X^*}(U)|q_\varphi(U)] = E[q_{X^*}(U)|\varphi],$$

since $F_\varphi(q_\varphi(t)-) \leq t \leq F_\varphi(q_\varphi(t))$ for all $t$. Lemma 2.3 and the law-invariance of $\rho$ imply that

$$\rho(-X) \leq \rho(-q_{X^*}(U)) = \rho(-X^*).$$

Moreover, the upper Hardy–Littlewood inequality and (21) yield that

$$\begin{aligned}
v \leq E[\varphi X^*] &\leq \int_0^1 q_\varphi(t) q_{X^*}(t) \, dt \\
&= E[q_\varphi(U) q_{X^*}(U)] = E[\varphi X],
\end{aligned} \tag{22}$$

and so $X$ solves the Neyman–Pearson problem, too. In view of Lemma 2.4, all inequalities in (22) must be identities. Hence, if $X^*$ is $\sigma(\varphi)$-measurable, then the "only if" part of Theorem 2.6 shows that $X^* = f^*(\varphi)$ for some increasing function $f^*$. $\square$

The argument in the preceding proof also yields the following uniqueness result for price densities with a continuous distribution. Remark 2.2 shows that this condition cannot be dropped.

PROPOSITION 2.7. *If $\varphi$ has a continuous distribution, then uniqueness in the class of $\sigma(\varphi)$-measurable solutions implies uniqueness in the class of $\mathcal{F}$-measurable solutions.*

PROOF. Let $X^*$ be an arbitrary solution and define

$$f := q_{X^*} \circ F_\varphi \quad \text{and} \quad X := f(\varphi).$$

Then $X$ has the same distribution as $X^*$. As in the proof of Proposition 2.5, we get that $X$ is a $\sigma(\varphi)$-measurable solution. Moreover, $E[X^*|\varphi]$ is also a $\sigma(\varphi)$-measurable solution by Lemma 2.3. Uniqueness gives $X = E[X^*|\varphi]$, and so $X^*$ has the same law as $E[X^*|\varphi]$. Hence,

$$0 = E[(X^*)^2] - E[E[X^*|\varphi]^2] = E[(X^* - E[X^*|\varphi])^2],$$

and we get that $P$-a.s. $X^* = E[X^*|\varphi]$. $\square$

Finally, we will need some elementary properties of the *minimal risk*

$$R_\varphi(v) := \min\{\rho(-X) \mid 0 \leq X \leq K, E[\varphi X] \geq v\}, \qquad 0 \leq v \leq K. \tag{23}$$

LEMMA 2.8. *$v \mapsto R_\varphi(v)$ is a continuous convex function that strictly increases from $\rho(0)$ to $\rho(-K)$ as $v$ increases from $0$ to $K$.*

PROOF. Clearly, $R_\varphi(0) = \rho(0)$ and $R_\varphi(K) = \rho(-K)$, due to our assumption $P[\varphi > 0] = 1$. It is also clear that $R_\varphi(v)$ is increasing in $v$. But if $R_\varphi(v) = R_\varphi(v')$ for some $v' > v$, then a solution for the Neyman–Pearson problem with $v'$ would also be a solution for $v$, a contradiction to Lemma 2.4. Therefore, the function $v \mapsto R_\varphi(v)$ is strictly increasing. Convexity easily follows from (11) and in turn implies continuity in the interior of $[0, K]$. Using (13), right-continuity at $v = 0$ follows from $R_\varphi(v) \le \rho(-v)$, while left-continuity at $v = K$ follows from (18). □

**3. Robust utility functionals defined in terms of density bounds.** For $\lambda \in (0, 1]$, let

$$\mathcal{Q}_\lambda = \left\{ Q \ll P \,\Big|\, \frac{dQ}{dP} \le \frac{1}{\lambda} P\text{-a.s.} \right\}$$

and note that $\mathcal{Q}_1 = \{P\}$. In this section, we solve the Neyman–Pearson problem for risk measures derived from robust utility functionals of the form

$$U_\lambda(X) = \min_{Q \in \mathcal{Q}_\lambda} E_Q[u(X)],$$

where $u$ is a utility function. Such utility functionals arise in a natural way from a robust Savage representation of preferences on asset profiles; see [14] and Section 2.5 in [13]. We will assume throughout this section that $u$ is concave, strictly increasing and continuously differentiable on its domain, which shall contain $[0, K]$. When measuring risk rather than utility, it is natural to switch signs and to introduce the convex increasing loss function $\ell(x) := -u(-x)$. Thus, we will consider the risk measure

(24) $$\rho_\lambda(-X) := -U_\lambda(-X) = \max_{Q \in \mathcal{Q}_\lambda} E_Q[\ell(X)].$$

If $\ell(x) = x$ for all $x$, then $\rho_\lambda$ reduces to the average value at risk $\text{AVaR}_\lambda$ of (8). The terminology "average value at risk" stems from the crucial fact that $\text{AVaR}_\lambda$ can be represented as an average of the upper values of the quantile function $q_X$ (the "value at risk") of $X \in L^\infty$:

(25) $$\text{AVaR}_\lambda(-X) = \max_{Q \in \mathcal{Q}_\lambda} E_Q[X] = \frac{1}{\lambda} \int_{1-\lambda}^1 q_X(t)\, dt;$$

see, for example, Theorem 4.39 in [13] and recall that we have assumed that $(\Omega, \mathcal{F}, P)$ is atomless. Thus, both $\text{AVaR}_\lambda$ and $\rho_\lambda$ are law-invariant and satisfy the general assumptions of Section 2.

We will first consider the Neyman–Pearson problem for the risk measure $\rho_\lambda$ of (24) in the case where the loss function $\ell$ is strictly convex on $[0, K]$; the Neyman–Pearson problem for $\text{AVaR}_\lambda$ will be considered in the next section. For simplicity, we will assume that the price density $\varphi$ is unbounded from above. Under our assumptions on the loss function $\ell$, its derivative $\ell'$ is a bijective function from its domain to some interval $(a,b) \subset (0,\infty)$. We extend its inverse function to all of $\mathbb{R}$ by setting

$$I(x) = \begin{cases} +\infty, & \text{for } x \geq b, \\ (\ell')^{-1}(x), & \text{for } a < x < b, \\ -\infty, & \text{for } x \leq a. \end{cases}$$

In the classical case $\lambda = 1$, we have $\mathcal{Q}_1 = \{P\}$, and it is well known that the unique solution of the Neyman–Pearson problem for $\rho_1$ takes the form

$$X_1^* = 0 \vee I(c_1 \varphi) \wedge K = (I(c_1(\varphi \vee y_1)) - I(c_1 y_1)) \wedge K,$$

where $c_1$ is the unique constant such that $E[\varphi X_1^*] = v$; see, for example, Section 3.3 of [13]. The parameter $y_1 = \ell'(0)/c_1$ can be interpreted as that level of prices at which the investor starts taking risky bets since the solution $X_1^*$ is supported on $\{\varphi \geq y_1\}$. Clearly, $c_1 = c_1(v)$ increases continuously from 0 to $+\infty$ and $y_1 = y_1(v)$ decreases continuously from $+\infty$ to 0 as $v$ increases from 0 to $K$.

For $0 < \lambda < 1$, we will see in the following theorem that $X_1^*$ also solves the Neyman–Pearson problem for $\rho_\lambda$, but only as long as the capital constraint $v$ does not exceed a certain *critical value* $v_\lambda$. For $v > v_\lambda$, a *diversification effect* will occur: the optimal solution will be a combination of a constant $\beta \in (0, v)$ and a classical solution $\widetilde{X}$ with upper bound $K - \beta$ and capital constraint $v - \beta$. Moreover, for all values of $v$, the classical part $\widetilde{X}$ will be concentrated on a subset of $\{\varphi \geq q\}$, where

$$q := q_\varphi(1 - \lambda).$$

Viewing the constant $\beta$ as a risk-free loan and $\widetilde{X}$ as a risky bet, we see that this effect is related to an "aversion" of the investor to accept risky bets on scenarios $\omega$ corresponding to prices $\varphi(\omega)$ which are not high enough. Hence, capital that cannot be raised by issuing a risky bet on high-price scenarios must instead be obtained via a risk-free loan. Note our shorthand notation of writing $E[X; \varphi \in A]$ for $E[X \mathbb{1}_{\{\varphi \in A\}}]$.

THEOREM 3.1. *Suppose that the distribution function of $\varphi$ is continuous and strictly increasing on $(0, \infty)$. Then*:

(a) *The Neyman–Pearson problem for $\rho_\lambda$ of* (24) *has a unique solution $X^*$ which is $P$-a.s. of the form*

(26) $$X^* = \beta + (I(c(\varphi \vee y)) - I(cy)) \wedge (K - \beta),$$

*where $\beta$, $y$ and $c$ are constants such that $\beta \geq 0$, $y \geq q$ and $c = \ell'(\beta)/y$.*

(b) *For every $\lambda \in (0,1)$, there exists a critical value $v_\lambda \in (0, KE[\varphi; \varphi \geq q])$ such that $\beta = 0$ if $v \leq v_\lambda$ and $0 < \beta < v$ for $v > v_\lambda$.*

(c) *The parameters $\beta$, $c$ and $y$ are increasing functions of $v$.*

PROOF. (a) By Proposition 2.5, we may concentrate on random variables $X$ that are of the form $X = f(\varphi)$ for an increasing function $f$. Then (20) and (25) imply that

$$
\begin{aligned}
\lambda \max_{Q \in \mathcal{Q}_\lambda} E_Q[\ell(X)] &= \int_{1-\lambda}^{1} q_{\ell(f(\varphi))}(t) \, dt \\
&= \int_{1-\lambda}^{1} \ell(f(q_\varphi(t))) \, dt \\
&= \int_{0}^{1} \ell(f(q_\varphi(t))) \mathbb{I}_{\{q_\varphi(t) \geq q\}} \, dt \\
&= E[\ell(f(\varphi)); \varphi \geq q],
\end{aligned}
\tag{27}
$$

where, in the third identity, we have used our assumptions on $\varphi$. Hence $X^* = f^*(\varphi)$ will solve the Neyman–Pearson problem provided that $f^*$ solves

(28) minimize $E[\ell(f(\varphi)); \varphi \geq q]$ among all increasing functions $f$ with $0 \leq f \leq K$ and $E[\varphi f(\varphi)] = v$, and vice versa. In particular, (28) admits a solution. It is clear that any such solution $f^*$ must satisfy $f^*(x) = f^*(q+)$ for all $x \leq q$. Taking $\beta := f^*(q)$ as given, the restriction of $f^*$ to $[q, \infty)$ is the unique solution to the following problem:

(29) minimize $E[\ell(f(\varphi)); \varphi \geq q]$ among all increasing functions $f$ on $[q, \infty)$ with $\beta \leq f \leq K$ and $E[\varphi f(\varphi); \varphi \geq q] = v - \beta E[\varphi; \varphi \leq q] =: v_\beta$.

If we drop the condition that $f$ in (29) is increasing, then it is well known (see, e.g., Section 3.3 of [13]) that (29) is solved by the function

$$f_\beta^*(x) = \beta \vee I(cx) \wedge K, \qquad x \geq q, \tag{30}$$

where $c$ is such that $E[\varphi f_\beta^*(\varphi); \varphi \geq q] = v_\beta$. But $f_\beta^*$ *is* increasing and hence solves (29). Since $\beta = f^*(q) = f^*(q+)$, we get $f^*(x) = \beta \vee I(cx) \wedge K$ for all $x \geq 0$. Moreover, there must be some $y \geq q$ such that $\beta = I(cy)$. Thus, $f^*$ can be written as $f^*(x) = \beta + (I(c(x \vee y)) - I(cy)) \wedge (K - \beta)$.

As for the uniqueness of solutions, we have just shown that all $\sigma(\varphi)$-measurable solutions are of the form (30) and can be parameterized via $\beta$. But a different $\beta$ needs a different $c$, so that two $\sigma(\varphi)$-measurable solutions must differ almost everywhere. The strict convexity of $\ell$ hence implies uniqueness of (28) and in turn uniqueness of the $\sigma(\varphi)$-measurable solution of the Neyman–Pearson problem. General uniqueness follows from Proposition 2.7.

Part (b) is obtained by combining Lemmas 3.2–3.4. Part (c) follows from Lemma 3.2 and the fact that $\beta < v$ as proved in Lemma 3.4. □

LEMMA 3.2. *The solutions in* (26) *are pointwise increasing in* $v$.

PROOF. Let $v$ and $v'$ be such that $0 \leq v < v' \leq K$, and consider the corresponding solutions $X^*(v)$ and $X^*(v')$. We want to show that $P$-a.s. $X^*(v') \geq X^*(v)$. To this end, define $X := X^*(v) \wedge X^*(v')$, $Y := X^*(v) - X$ and $Z := X^*(v') - X$. Then $v_0 := E[\varphi X] \leq v$, and there exists $\alpha \in (0,1]$ such that $(1-\alpha)E[\varphi Z] = E[\varphi Y] = v - v_0$. Clearly, we have $Y = 0$ on $\{Z > 0\}$ and hence, by the convexity of $\ell$, $P$-a.s.,

$$\begin{aligned}
\ell(X + (1-\alpha)Z + \alpha Z) - \ell(X + Y + \alpha Z) \\
\geq \ell(X + (1-\alpha)Z) - \ell(X + Y).
\end{aligned} \tag{31}$$

Both $X^*(v)$ and $X^*(v')$ are increasing functions of the price density $\varphi$, and one easily checks that the same is true of $X + Y + \alpha Z$ and of $X + (1-\alpha)Z$. Hence, multiplying (31) with $\mathbb{I}_{\{\varphi \geq q\}}$, taking expectations with respect to $P$, and using (27) yields

$$\begin{aligned}
\max_{Q \in \mathcal{Q}_\lambda} E_Q[\ell(X^*(v'))] &- \max_{Q \in \mathcal{Q}_\lambda} E_Q[\ell(X^*(v) + \alpha Z)] \\
&\geq \max_{Q \in \mathcal{Q}_\lambda} E_Q[\ell(X + (1-\alpha)Z)] - \max_{Q \in \mathcal{Q}_\lambda} E_Q[\ell(X^*(v))] \\
&\geq 0,
\end{aligned} \tag{32}$$

where the latter inequality follows from the fact that $E[\varphi(X + (1-\alpha)Z)] = v$. Moreover, $E[\varphi(X^*(v) + \alpha Z)] = v'$, which in view of (32) and the uniqueness of solutions implies that $P$-a.s. $X^*(v') = X^*(v) + \alpha Z \geq X^*(v)$. □

LEMMA 3.3. *For every $\lambda \in (0,1)$, there exists $\varepsilon > 0$ such that $\beta = 0$ for $v \leq \varepsilon$.*

PROOF. Fix $v \in (0, K)$ for the first step. For $\gamma \in [0, v)$, let $f_\gamma(x) := \gamma \vee I(c_\gamma x) \wedge K$, where $c_\gamma \in (0, \infty)$ is such that $E[\varphi f_\gamma(\varphi)] = v$. We denote by $y^\gamma := \ell'(\gamma)/c_\gamma$ the point at which $f_\gamma$ starts being larger than $\gamma$. Suppose that $\gamma' > \gamma$. Then $y^1 := y^\gamma \wedge y^{\gamma'} > 0$ and

$$\begin{aligned}
E[\varphi \cdot (\gamma' \vee I(c_{\gamma'}\varphi) \wedge K); \varphi > y^1] \\
= v - \gamma' E[\varphi; \varphi \leq y^1] \\
< E[\varphi f_\gamma(\varphi); \varphi > y^1] \\
\leq E[\varphi \cdot (\gamma' \vee I(c_\gamma \varphi) \wedge K); \varphi > y^1].
\end{aligned}$$

It follows that $\gamma \mapsto c_\gamma = c_\gamma(v)$ is strictly decreasing and that $\gamma \mapsto y^\gamma = y^\gamma(v)$ is strictly increasing as long as $v$ is fixed.

Now let $L(\gamma) := E[\ell(f_\gamma(\varphi)); \varphi \geq q]$. It follows from the proof of Theorem 3.1 that $\beta = 0$ or $\beta = v$ if

(33) $$L(\gamma) - L(0) > 0 \qquad \text{for all } \gamma \in (0, v).$$

But (33) also implies that $L(v) := \ell(v)P[\varphi \geq q] = \lim_{\gamma \uparrow v} L(\gamma) > L(0)$, for the case $L(v) = L(0)$ is excluded by the uniqueness of the solution (26). Hence, (33) is equivalent to $\beta = 0$.

In addition to $y^\gamma$, we will also need the point $y_\gamma := \ell'(\gamma)/c_0 < y^\gamma$ at which $f_0$ leaves the level $\gamma$. Letting $\Delta := \ell(f_\gamma(\varphi)) - \ell(f_0(\varphi))$, we have

$$L(\gamma) - L(0) = E[\Delta; \varphi \geq y_\gamma] + E[\Delta; y_0 \leq \varphi < y_\gamma]$$
$$+ E[\Delta; q \leq \varphi < y_0].$$

On $\{\varphi \geq y_\gamma\}$, we get from the first step that $f_\gamma(\varphi) \leq f_0(\varphi)$ and in turn

$$\Delta \geq \ell'(f_0(\varphi))[f_\gamma(\varphi) - f_0(\varphi)] \geq c_0\varphi[f_\gamma(\varphi) - f_0(\varphi)].$$

Moreover, $\Delta \geq 0$ on $\{y_0 \leq \varphi < y_\gamma\}$, and on $\{q \leq \varphi < y_0\}$ we have $\gamma = f_\gamma(\varphi) \geq f_0(\varphi) = 0$. Therefore,

$$L(\gamma) - L(0) \geq c_0 E[\varphi \cdot (f_\gamma(\varphi) - f_0(\varphi)); \varphi \geq y_\gamma]$$
$$+ (\ell(\gamma) - \ell(0))P[q \leq \varphi < y_0]$$
$$\geq c_0(v - \gamma P[\varphi < y_\gamma] - v) + \gamma \ell'(0)P[q \leq \varphi < y_0]$$
$$\geq \gamma(\ell'(0)P[q \leq \varphi < \ell'(0)/c_0] - c_0).$$

By our assumption that $\varphi$ has a continuous and strictly increasing distribution function, the factor $c_0 = c_0(v)$ tends continuously from 0 to $+\infty$ as $v$ increases from 0 to $K$, and so the right-hand side will be strictly positive as soon as $v$ is small enough and $\gamma$ is between 0 and $v$. $\square$

LEMMA 3.4. *We have $\beta < v$ for all $v \in (0, K)$ and $\beta > 0$ for $v > KE[\varphi; \varphi \geq q]$.*

PROOF. As to the first part of the assertion, it follows from Lemma 3.3 that $X \equiv v$ is not optimal for small enough $v > 0$. That is, $R_\varphi(v) < \rho(-X) = v$, where $R_\varphi(v)$ is as in (23). The convexity of $v \mapsto R_\varphi(v) - v$, which follows from Lemma 2.8, hence implies that $v = 0$ and $v = K$ are the only two points in $[0, K]$ with $R_\varphi(v) = v$. Thus, $X \equiv v$ cannot be optimal for any $v \in (0, K)$.

The second part of the assertion follows immediately from the fact that the parameter $y$ in Theorem 3.1 has been shown to be larger than or equal to $q$. $\square$

Let us now briefly comment on the translation invariant modification

$$\widehat{\rho}_\lambda(-X) = \inf\left\{m \in \mathbb{R} \,\bigg|\, \max_{Q \in \mathcal{Q}_\lambda} E_Q[\ell(X - m)] \leq x_0\right\}$$

of $\rho_\lambda$, which is a convex measure of risk in the sense of [11]. In addition to the assumptions made at the beginning of this section, we assume that $\ell$ is defined on all of $\mathbb{R}$, and $x_0$ is a fixed interior point of $\ell(\mathbb{R})$. Clearly, $\widehat{\rho}_\lambda$ is law-invariant and satisfies the properties (10) through (16). We denote by $R(v) := R_\varphi(v)$ the minimal risk for $\widehat{\rho}_\lambda$, as defined in (23). Recall that $q$ denotes the $(1-\lambda)$-quantile of $\varphi$.

COROLLARY 3.5. *Suppose that the distribution function of $\varphi$ is continuous and strictly increasing on $(0, \infty)$. Then the Neyman–Pearson problem for $\widehat{\rho}_\lambda$ has a unique solution $X^*$ that is $P$-a.s. of the form*

$$X^* = \alpha + (I(\gamma(\varphi \vee z)) - I(\gamma z)) \wedge (K - \alpha),$$

*where $\alpha$, $z$ and $\gamma$ are constants such that $0 \leq \alpha < v$, $z \geq q$ and $\gamma = \ell'(\alpha - R(v))/z$. Moreover, for every $\lambda \in (0, 1)$, there exists a critical value $\widehat{v}_\lambda \in (0, K)$ such that $\alpha = 0$ if $v \leq \widehat{v}_\lambda$.*

PROOF. Take a solution $X^*$ at level $v$ and let $\ell_{R(v)}(x) := \ell(x - R(v))$. Then we see that $\max_{Q \in \mathcal{Q}_\lambda} E[\ell_{R(v)}(X^*)] = x_0$. On the other hand, if $0 \leq X \leq K$ and $E[\varphi X] \geq v$ but $X$ is *not* a solution, then we must have $\widehat{\rho}_\lambda(-X) > R(v)$ and hence $\max_{Q \in \mathcal{Q}_\lambda} E[\ell_{R(v)}(X)] > x_0$. So $X^*$ solves the Neyman–Pearson problem for $\widehat{\rho}_\lambda$ at level $v$ if and only if $X^*$ minimizes $\max_{Q \in \mathcal{Q}_\lambda} E[\ell_{R(v)}(X)]$ among all $X$ with $0 \leq X \leq K$ and $E[\varphi X] \geq v$. For fixed $v$, this problem is of the same type as the one of Theorem 3.1, and so we get a representation of solutions in terms of the inverse $I_{R(v)}$ of $\ell'_{R(v)}$. But $I_{R(v)}(x) = I(x) + R(v)$, and we obtain the first part of the assertion. The existence of the critical value $\widehat{v}_\lambda$ follows by the same arguments as in Lemma 3.3 when one replaces $\ell$ by $\ell'_{R(v)}$ and $I$ by $I_{R(v)}$; only minor modifications are needed. □

From the proof it is clear that, for given $v > 0$, the parameters $\alpha$, $\gamma$ and $z$ will generally be different from the corresponding parameters $\beta$, $c$ and $y$ in Theorem 3.1, because the problem now involves the loss function $\ell_{R(v)}(x) := \ell(x - R(v))$ rather than $\ell$ itself. Also, in the case in which $\alpha = 0$ but $\lambda < 1$, the solution $X^*_\lambda := X^*$ typically does *not* coincide with the solution $X^*_1$ to the Neyman–Pearson problem for the "classical" risk measure

$$\widehat{\rho}_1(-X) = \inf\{m \in \mathbb{R} \mid E[\ell(X - m)] \leq x_0\}$$

[with the exception of an exponential loss function $\ell(x) = e^{\alpha x}$]. To see this, note first that

$$\max_{Q \in \mathcal{Q}_\lambda} E_Q[\ell(X - m)] > E[\ell(X - m)]$$

unless $X$ is constant. This in turn implies that $R(v) = \widehat{\rho}_\lambda(-X_\lambda^*) > \widehat{\rho}_1(-X_1^*) =: R^1(v)$ for otherwise $\widehat{\rho}_1(-X_\lambda^*)$ would be strictly less than $\widehat{\rho}_1(-X_1^*)$. But $X_\lambda^*$ is of the form

$$X_\lambda^* = 0 \vee (I(\gamma_\lambda \varphi) + R(v)) \wedge K,$$

while

$$X_1^* = 0 \vee (I(\gamma_1 \varphi) + R^1(v)) \wedge K,$$

which shows that $\gamma_\lambda < \gamma_1$.

**4. Quantile-based coherent risk measures.** A quantile-based coherent risk measure is of the form

$$\rho_k(-X) := \int_0^1 k(t) q_X(t)\, dt, \qquad X \in L^\infty,$$

where $k : [0,1) \to [0, \infty)$ is an increasing right-continuous function such that $\int_0^1 k(t)\, dt = 1$, and where $q_Y$ denotes a quantile function of a random variable $Y$. The average value at risk $\mathrm{AVaR}_\lambda$ of (25) is thus the particular quantile-based coherent risk measure with $k = \frac{1}{\lambda} \mathbb{I}_{[1-\lambda, 1)}$. For general $k$, let $\widetilde{\mu}$ be the positive Radon measure on $[0,1)$ such that $k(t) = \widetilde{\mu}([0,t])$. Then $\mu(d\lambda) = (1-\lambda)\widetilde{\mu}(d\lambda)$ is a probability measure on $[0,1)$ such that

$$\rho_k(-X) = \int_{[0,1)} \mathrm{AVaR}_{1-\lambda}(-X) \mu(d\lambda).$$

Since $\mathrm{AVaR}_\lambda$ is a coherent measure of risk which is continuous from below and, hence, from above (see, e.g., [20] or Theorem 4.39 in [13]), the same is true of the quantile-based risk measure $\rho_k$. In particular, $\rho_k$ satisfies the properties (10) through (16) and can be represented in the form

$$\rho_k(-X) = \max_{Q \in \mathcal{Q}^k} E_Q[X],$$

where $\mathcal{Q}^k$ is a set of probability measures, which has been described by Dana and Carlier [5].

Let us now turn to the Neyman–Pearson problem for $\rho_k$. By the positive homogeneity of $\rho_k$, there is no loss in generality if we assume that $K = 1$. Our first result in this section will show that the Neyman–Pearson problem for $\rho_k$ can be reduced to the minimization of an ordinary expectation over a very limited class $\mathcal{J}$ of functions. This class $\mathcal{J}$ consists of all increasing step functions $f : (0, \infty) \to [0, 1]$ that take at most one value in $(0, 1)$. More precisely, each $f \in \mathcal{J}$ can be written as

$$f = \beta \mathbb{I}_{J_0} + \mathbb{I}_{J_1}$$

for some $\beta \in (0, 1)$ and two disjoint intervals $J_0, J_1 \subset (0, \infty)$ such that $J_0$ is either empty or satisfies $P[\varphi \in J_0] > 0$. Here and in the sequel, we use the term interval in a broad sense: an interval may also be empty or consisting of a single element. Since $f$ must be increasing, $J_1$ must either be empty or unbounded to the right. If both $J_0$ and $J_1$ are nonempty, then the right-hand endpoint of $J_0$ must coincide with the left-hand endpoint of $J_1$.

Recall that $F_\varphi$ denotes the distribution function of $\varphi$ under $P$ and let us introduce the function

$$g_k(x) = \begin{cases} k(F_\varphi(x)), & \text{if } F_\varphi \text{ is continuous at } x, \\ \dfrac{1}{F_\varphi(x) - F_\varphi(x-)} \int_{F_\varphi(x-)}^{F_\varphi(x)} k(t)\, dt, & \text{otherwise.} \end{cases}$$

Consider the following variational problem:

(34) $\quad$ minimize $E[g_k(\varphi) f(\varphi)]$ among all increasing functions $f$ with $0 \leq f \leq 1$ and $E[\varphi f(\varphi)] = v$.

It would be tempting to apply the classical Neyman–Pearson lemma to solving (34), but this approach would only work if the function $g_k(x)/x$ were *decreasing* in $x$, because otherwise we might not obtain an *increasing* solution $f$.

THEOREM 4.1.

(a) If $f^*$ solves (34), then $X^* := f^*(\varphi)$ solves the Neyman–Pearson problem for $\rho_k$.

(b) There exists a function $f^* \in \mathcal{J}$ that solves (34).

(c) If $f^*$ is such that $f^*(\varphi)$ solves the Neyman–Pearson problem for $\rho_k$, then $f^*$ solves (34).

(d) If the solution $f^* \in \mathcal{J}$ of part (b) is unique within $\mathcal{J}$ up to $(P \circ \varphi^{-1})$-nullsets, then $X^* = f^*(\varphi)$ is the $P$-a.s. unique $\sigma(\varphi)$-measurable solution to the Neyman–Pearson problem for $\rho_k$.

The proof of Theorem 4.1 is deferred to Section 5. Here we will first illustrate how this result leads to explicit solutions of the Neyman–Pearson problem for quantile-based coherent risk measures. In order not to complicate the presentation, we assume for the rest of this section that the distribution function $F_\varphi$ is continuous and strictly increasing on $\{x > 0 \mid F_\varphi(x) < 1\}$. Then the corresponding quantile function $q_\varphi$ will also be continuous and strictly increasing. We let $q_\varphi(0) := 0$ and $q_\varphi(1) := \|\varphi\|_{L^\infty} \leq \infty$, and we define two functions $\Phi$ and $\Gamma$ by

$$\Phi(x) := \int_0^x q_\varphi(t)\, dt \quad \text{and} \quad \Gamma(x) := \int_0^x k(t)\, dt, \qquad 0 \leq x \leq 1.$$

Then we take the unique $z_v$ such that

$$\Phi(z_v) = 1 - v,$$

and define two functions $\beta$ and $R$ on $\Delta_v := \{(x,y) \mid 0 \leq x < z_v < y \leq 1\} \cup \{(z_v, z_v)\}$ by

$$\beta(x,y) := \begin{cases} 0, & \text{if } x = z_v = y, \\ \dfrac{v - 1 + \Phi(y)}{\Phi(y) - \Phi(x)}, & \text{otherwise,} \end{cases}$$

$$R(x,y) := \beta(x,y)[\Gamma(y) - \Gamma(x)] + 1 - \Gamma(y).$$

COROLLARY 4.2. *Suppose that the pair $(x^*, y^*)$ minimizes the function $R$ over the domain $\Delta_v$, and let $\beta^* := \beta(x^*, y^*)$, $a := q_\varphi(x^*)$ and $b := q_\varphi(y^*)$. Then $X^* := f^*(\varphi)$ solves the Neyman–Pearson problem for $\rho_k$, where*

(35) $$f^* := \beta^* \mathbb{I}_{[a,b)} + \mathbb{I}_{[b,\infty)}.$$

*Conversely, suppose that $f \in \mathcal{J}$ is a.e. of the form (35) and solves (34). Then the pair $(x^*, y^*) := (F_\varphi(a), F_\varphi(b))$ minimizes $R$ on $\Delta_v$. In particular, the Neyman–Pearson problem for $\rho_k$ has a unique solution if and only if $R$ has a unique minimizer on $\Delta_v$.*

PROOF. It is straightforward to verify that a function $f = \beta \mathbb{I}_{[a,b)} + \mathbb{I}_{[b,\infty)} \in \mathcal{J}$ satisfies the constraints in (34) if and only if $(x,y) := (F_\varphi(a), F_\varphi(b)) \in \Delta_v$ and $\beta = \beta(x,y)$. An analogous computation shows that $E[g_k(\varphi) f(\varphi)] = R(x,y)$, so that the assertion follows from Theorem 4.1 and Proposition 2.7. □

The preceding corollary implies that $\sigma(\varphi)$-measurable solutions to the Neyman–Pearson problem need not be unique, even for genuinely nonadditive risk measures and for price densities with a continuous distribution.

REMARK 4.3. In the case $k \equiv q_\varphi$, we have $R(x,y) = v$ for all $(x,y) \in \Delta_v$. Hence, *each* function

$$f = \beta(x,y)\mathbb{I}_{[a,b)} + \mathbb{I}_{[b,\infty)} \qquad \text{for } a = q_\varphi(x), \ b = q_\varphi(y)$$

solves the Neyman–Pearson problem for $\rho_k$, and so does every convex combination of these functions.

Below, we will use Corollary 4.2 to obtain an explicit solution for the Neyman–Pearson problem for $\text{AVaR}_\lambda$. As one may guess from Theorem 3.1, we will find the dichotomy $x^* = y^* = z_v$ or $x^* = 0$ and $y^* > z_v$. But before doing so, let us show in the following example that the case $0 < x^* < y^* < 1$ can also occur.

EXAMPLE 4.4. Let us consider the case in which $\varphi$ has a uniform distribution on $(0,2)$, so that $q_\varphi(t) = 2t$, $\Phi(x) = x^2$ and $z_v = \sqrt{1-v}$. We take

$$k = \tfrac{1}{2}\mathbb{I}_{[0,\xi)} + \lambda \mathbb{I}_{[\xi,1)},$$

where $\xi \in (\tfrac{1}{2}, 1)$ and $\lambda$ is such that $k$ integrates to 1. With this choice, $\Gamma(x) < x^2 = \Phi(x)$ for all $x \in (\tfrac{1}{2}, \xi]$. Consequently, $R(z_v, z_v) = 1 - \Gamma(z_v) > v = R(0,1)$ for all $1 - \xi^2 \leq v < 3/4$. It follows that $(z_v, z_v)$ does not minimize $R$ for those values of $v$. Let $(x^*, y^*)$ be a minimizer of $R$ on $\Delta_v$. Then the right-hand derivative of $x \mapsto R(x, y^*)$ is equal to

$$\beta(x, y^*)\left(2x \frac{\Gamma(y^*) - \Gamma(x)}{(y^*)^2 - x^2} - k(x)\right).$$

Since this expression is strictly negative for small enough $x$, the optimal $x^*$ must be larger than 0.

Let us now show that the case $y^* = 1$ cannot occur if the parameter $\xi$ is sufficiently close to $\tfrac{1}{2}$ and $1 - \xi^2 \leq v < \tfrac{3}{4}$. To this end, one verifies first that the left-hand derivative of $y \mapsto R(x^*, y)$ at $y = 1$ is given by

(36) $$(1 - \beta(x^*, 1))\left(2\frac{1 - \Gamma(x^*)}{1 - (x^*)^2} - \lambda\right).$$

For $0 \leq x \leq z_v$ and $z_v \leq \xi$, the function $(1 - \Gamma(x))/(1 - x^2)$ has a global minimum at $x = \tfrac{1}{4}$, where it takes the value $\tfrac{14}{15}$. On the other hand, $\lambda$ tends to $\tfrac{3}{2}$ when $\xi$ goes to $\tfrac{1}{2}$. Thus, (36) must be strictly positive if $\xi$ is not too large, and we conclude that $y = 1$ cannot be optimal.

Let us now turn to the Neyman–Pearson problem for $\mathrm{AVaR}_\lambda$. There are various ways of handling this special case. For instance, one can use the arguments of the proof of Theorem 3.1 to reduce the problem to the variational problem (29) for $\ell(x) = x$, which can then be solved via the classical Neyman–Pearson lemma. Here we will use instead a computation based on Corollary 4.2.

As in Theorem 3.1, we will find a critical value $v_\lambda$ such that the solution reduces to the solution for $\rho(-X) = E[X]$ as long as $v \leq v_\lambda$. That is, the solution provided by the classical Neyman–Pearson lemma is optimal for capital levels $v \leq v_\lambda$. For $v > v_\lambda$, the solution will be a nontrivial convex combination of the classical solution at level $v_\lambda$ and of a risk-free unit investment. This critical value will be of the form

$$v_\lambda = 1 - \Phi(y_\lambda),$$

where $y_\lambda \in (1 - \lambda, 1]$ is defined as the unique maximizer of the function

$$(0,1] \ni y \longmapsto \frac{y + \lambda - 1}{\Phi(y)}.$$

Thus, if $q_\varphi(1) = \|\varphi\|_{L^\infty} > \lambda^{-1}$, then $y_\lambda \in (1-\lambda, 1)$ is the unique solution to the equation

$$q_\varphi(y_\lambda)(y_\lambda + \lambda - 1) = \Phi(y_\lambda).$$

COROLLARY 4.5. *The Neyman–Pearson problem for* $\mathrm{AVaR}_\lambda$ *has a unique solution* $X^*$. *If* $v \leq v_\lambda$, *then*

$$X^* = \mathbb{I}_{[b_0, \infty)}(\varphi),$$

*where* $b_0 := q_\varphi(z_v)$. *If* $v > v_\lambda$, *then the solution is given by*

$$X^* = \beta^* + (1-\beta^*)\mathbb{I}_{[b_1, \infty)}(\varphi),$$

*where* $\beta^* = \beta(0, y_\lambda)$ *and* $b_1 = q_\varphi(y_\lambda)$. *Moreover, with* $C_\lambda := (y_\lambda + \lambda - 1)/\Phi(y_\lambda)$, *the minimal risk* (23) *is given by*

$$R_\varphi(v) = \begin{cases} (1-z_v)/\lambda, & \text{if } v \leq v_\lambda, \\ 1 - C_\lambda(1-v)/\lambda, & \text{if } v > v_\lambda. \end{cases}$$

PROOF. It suffices to consider the case $0 < v < 1$. For $k = \frac{1}{\lambda}\mathbb{I}_{[1-\lambda, 1)}$, we have

$$\lambda R(x, y) = \beta(x, y)[(y+\lambda-1) \vee 0 - (x+\lambda-1) \vee 0] + \lambda - (y+\lambda-1) \vee 0.$$

Let $(x^*, y^*)$ be a minimizer of $R$ on $\Delta_v$, and suppose first that $(x^*, y^*) \neq (z_v, z_v)$. Then $y^* > z_v$ and

$$\frac{\partial}{\partial x}\beta(x, y^*) = \beta(x, y^*)\frac{q_\varphi(x)}{\Phi(y^*) - \Phi(x)} > 0.$$

Thus, we see that $\beta(x, y^*)$ and, hence, $R(x, y^*)$ are strictly increasing in $x$ as long as $x < 1 - \lambda$. If, on the other hand, $x > 1 - \lambda$, then

$$\lambda \frac{\partial}{\partial x} R(x, y^*) = \beta(x, y^*)\left(\frac{q_\varphi(x)(y^* - x)}{\Phi(y^*) - \Phi(x)} - 1\right) < 0.$$

So $x^*$ must be equal to either 0 or $z_v$.

Let us now look for the optimal $y^*$, given that $x^* = 0$. We have

$$\lambda R(0, y) = \lambda - (1-v)\frac{(y+\lambda-1) \vee 0}{\Phi(y)}.$$

For $y \leq 1 - \lambda$, this yields $R(0, y) = 1$, which according to Lemma 2.4 cannot be optimal. For $y > 1 - \lambda$, the choice $(x^*, y^*) = (0, y_\lambda)$ will be optimal—but only if $y_\lambda > z_v$ and unless the alternative choice $(x^*, y^*) = (z_v, z_v)$ gives a better result.

If $y_\lambda \leq z_v$, then $y \mapsto R(0,y)$ has no minimizer on $(z_v, 1]$, and it follows that $(x^*, y^*) = (z_v, z_v)$ must be the optimal choice. Note that $y_\lambda > z_v$ if and only if $v > v_\lambda$.

Finally, let us compare $R(0, y_\lambda)$ against $R(z_v, z_v)$ in case that $y_\lambda > z_v$. Since $y_\lambda > 1 - \lambda$, we have

$$\lambda R(0, y_\lambda) = \lambda - (1-v)\frac{y_\lambda + \lambda - 1}{\Phi(y_\lambda)} = \lambda - \Phi(z_v)\frac{y_\lambda + \lambda - 1}{\Phi(y_\lambda)}$$

and

$$\lambda R(z_v, z_v) = \lambda - (z_v + \lambda - 1) \vee 0.$$

Since $y_\lambda$ is the unique maximizer of the function $x \mapsto (x + \lambda - 1)/\Phi(x)$, we thus see that $R(0, y_\lambda)$ is strictly better than $R(z_v, z_v)$ and hence $(x^*, y^*) = (0, y_\lambda)$ as long as $y_\lambda > z_v$. An application of Corollary 4.2 concludes the proof. □

REMARK 4.6 (Comparison with value at risk). Consider the value at risk at level $\lambda \in (0,1)$,

$$\text{VaR}_\lambda(-X) = \inf\{m \in \mathbb{R} \mid P[X > m] \leq \lambda\},$$

which is a quantile-based risk measure that satisfies all the assumptions of Section 2 except for convexity (11). Denoting $R_\varphi(v)$ the corresponding minimax risk (23), we see that $X^*$ solves the Neyman–Pearson problem for $\text{VaR}_\lambda$ if $P[X^* > R_\varphi(v)] \leq \lambda$ and $E[\varphi X^*] \geq v$. Thus, for $v$ with $z_v > q := q_\varphi(1 - \lambda)$, any $X$ that is concentrated on $\{\varphi > q\}$ and satisfies $0 \leq X \leq 1$ and $E[\varphi X] \geq v$ solves our problem and has risk $R_\varphi(v) = \text{VaR}_\lambda(-X) = 0$. For $z_v \leq q$, there is a unique solution of the form

$$X^* = r\mathbb{I}_{[0,q)}(\varphi) + \mathbb{I}_{[q,\infty)}(\varphi),$$

where $r = R_\varphi(v)$ is determined by the budget constraint $E[\varphi X^*] = v$. This solution is similar to the one for $\text{AVaR}_\lambda$, but involves different parameters.

REMARK 4.7. It follows from the results of Kusuoka [20] and Delbaen [7] that, for a quantile-based coherent risk measure $\rho_k$, the set function $v_k(A) := \rho_k(-\mathbb{I}_A)$ is a 2-alternating Choquet capaticity. Therefore, the Neyman–Pearson problem for $\rho_k$ falls within the range of the Neyman–Pearson lemma for capacities as proved by Huber and Strassen [16], and our results can be interpreted in terms of the Radon–Nykodym derivative $\pi$ of the measure $dP^* := \varphi \, dP$ with respect to the capacity $v_k$. In the case of $\text{AVaR}_\lambda$, we get for $\|\varphi\|_{L^\infty} > \lambda^{-1}$ that $\pi = c \cdot \varphi \vee q_\varphi(y_\lambda)$ for some constant $c > 0$. It is shown in [16] that $\pi = d\widehat{P}/dQ_0$ for some $Q_0 \in \mathcal{Q}_\lambda$, and we get $c \cdot \varphi \vee q_\varphi(y_\lambda) = \varphi \cdot dP/dQ$. Using our formulae for $y_\lambda$, one easily obtains $c = \lambda$, that is,

$$\pi = \lambda(\varphi \vee q_\varphi(y_\lambda)).$$

This extends earlier results by Rieder [23] and Bednarski [2].

## 5. Proof of Theorem 4.1.

PROOF OF PARTS (a) AND (c). As in the proof of Proposition 2.5, we see that $g_k(q_\varphi) = E_\lambda[k|q_\varphi]$, where $\lambda$ denotes the Lebesgue measure on $(0,1)$. Hence, for any increasing function $f:[0,\infty) \to [0,1]$,

$$\rho(-f(\varphi)) = \int_0^1 k(t) q_{f(\varphi)}(t)\, dt$$
$$= \int_0^1 g_k(q_\varphi(t)) f(q_\varphi(t))\, dt$$
$$= E[g_k(\varphi) f(\varphi)],$$

where we have used (20) in the second step. Applying Lemma 2.4 and Proposition 2.5 yields (a) and (c). $\square$

The proof of parts (b) and (d) requires some preparation. To illustrate our idea of solving (34), suppose first that the price density has a continuous distribution. Then we may, without loss of generality, restrict our attention to right-continuous increasing functions $f:[0,\infty) \to [0,1]$ in (34). Via $f(x) = \nu([0,x])$, any such function $f$ can be identified with a unique subprobability measure $\nu$ on $[0,\infty)$ and vice versa. Fubini's theorem implies that

$$E[g_k(\varphi) f(\varphi)] = \int G_k(x) \nu(dx)$$

and

$$E[\varphi f(\varphi)] = \int G_\varphi(x) \nu(dx),$$

where $G_k(x) = E[g_k(\varphi); \varphi \geq x]$ and $G_\varphi(x) = E[\varphi; \varphi \geq x]$. Thus, (34) is equivalent to minimizing the integral $\int G_k(x) \nu(dx)$ over the convex set $\widetilde{C}$ of all subprobability measures $\nu$ on $[0,\infty)$ that satisfy the constraint

$$\int G_\varphi(x) \nu(dx) = v.$$

Our strategy of solving this moment problem is to identify the extreme points of $\widetilde{C}$ as those subprobability measures that correspond to functions in $C \cap \mathcal{J}$ and to show that it suffices to minimize $\int G_k(x) \nu(dx)$ among such extreme measures $\nu$.

If the distribution of $\varphi$ is not continuous, the problem becomes slightly more involved. This is mainly due to the fact that we may no longer pass to a right-continuous (or left-continuous) version of $f$. We may only suppose that $f$ is right-continuous on the set $[0,\infty) \setminus D$, where $D$ denotes the set of discontinuity points of $F_\varphi$. Nevertheless, it will be possible to identify $f$ with a certain measure $\nu$ living on a larger space $S \supset [0,\infty)$, in which each point in $D$ occurs twice. Our problem (34) will then turn out to be equivalent to a certain moment problem for these measures $\nu$. Once this identification has been achieved, the extreme points of the set defined by the moment constraint on $\nu$ can be identified by using general results like those proved by Winkler [25]. In our simple situation, however, we will avoid using the general theory. Instead, we will give a short and straightforward argument in identifying the extreme points.

Define a probability measure $\mu$ on $[0, \infty)$ by

$$\mu(A) := E[\varphi; \varphi \in A],$$

and denote by $C$ the convex set of all increasing functions $f : [0, \infty) \to [0, 1]$ that are right-continuous on $D^c$ and satisfy the constraint $\int f \, d\mu = v$.

LEMMA 5.1. *The set of extreme points of $C$ is given by*

(37) $$\operatorname{ext} C = C \cap \mathcal{J}.$$

PROOF. First we show the inclusion $\supset$ in (37). So suppose that $f \in C \cap \mathcal{J}_\mu$ is of the form $f = \beta \mathbb{I}_{J_0} + \mathbb{I}_{J_1}$ and can be written as $f = \lambda f_1 + (1 - \lambda) f_2$ for certain $f_i \in C$ and $\lambda \in (0, 1)$. Since $0 \le f_i \le 1$, we get immediately that $f_i = 0$ on $(J_0 \cup J_1)^c$ and $f_i = 1$ on $J_1$. This proves that $f$ is an extreme point if $J_0$ is empty. Now suppose that $J_0$ is nonempty. Then $\mu(J_0) > 0$ by the definition of $\mathcal{J}$. Since both $f_1$ and $f_2$ are increasing, each $f_i$ must be equal to some constant $\beta_i \in [0, 1]$ on $J_0$. But then the conditions $\int f_i \, d\mu = v$ and $\mu(J_0) > 0$ imply $\beta_1 = \beta_2 = \beta$ and in turn $f_1 = f_2 = f$.

For the proof of the inclusion $\subset$ in (37), it will be convenient to identify a function $f \in C$ with a suitable subprobability measure $\nu$. To this end, we first define a subprobability measure $\nu_D$ on $D$ by

$$\nu_D := \sum_{x \in D} (f(x+) - f(x)) \delta_x.$$

Then we let

$$f_D(x) := \nu_D([0, x)) \quad \text{and} \quad f_c(x) := f(x) - f_D(x), \qquad x \ge 0.$$

Note that $f_c$ is right-continuous and increasing. Hence, there exists a subprobability measure $\nu_c$ on $[0, \infty)$ such that $f_c(x) = \nu_c([0, x])$.

Now think of the set $D$ as being *separate* from $[0, \infty)$, and consider the set $S := [0, \infty) \cup D$, on which every discontinuity point of $\mu$ is represented *twice*. Our function $f$ gives rise to a subprobability measure $\nu$ on $S$ defined for Borel sets $A \subset S$ by

$$\nu(A) := \nu_c(A \cap [0, \infty)) + \nu_D(A \cap D).$$

Conversely, any subprobability measure $\tilde{\nu}$ on the Borel field of $S$ gives rise to an increasing function $\tilde{f}$ on $[0, \infty)$ that is right-continuous except possibly at discontinuity points of $\mu$: simply let $\tilde f(x) := \tilde\nu(A_x)$, where $A_x := [0,x] \cup \{y \in D \,|\, y < x\}$. Note also that $f = \tilde f$ if and only if $\nu = \tilde\nu$.

By means of Fubini's theorem, we find that $\int f\,d\mu = \int G\,d\nu$, where $G$ is the function on $S$ defined by

$$\text{(38)} \qquad G(x) = \begin{cases} \mu([x,\infty)), & \text{for } x \in [0,\infty), \\ \mu((x,\infty)), & \text{for } x \in D. \end{cases}$$

Hence, $C$ can be identified with the set of all subprobability measures $\nu$ on $S$ such that $\int G\,d\nu = v$.

Let us now consider the case in which $\sup_x f(x) = 1$, corresponding to $\nu(S) = 1$. Suppose $f$ takes more than one value in $(0,1)$. Then $S$ can be decomposed into three disjoint sets $A_1, A_2, A_3$ such that $a_i := \nu(A_i) > 0$. Let also $b_i := \int_{A_i} G\,d\nu$, and denote by $\nu|_{A_i}$ the measure $\nu|_{A_i}(A) := \nu(A \cap A_i)$. For coefficients $\alpha_i \geq 0$, the measure

$$\nu_1 := \alpha_1 \nu|_{A_1} + \alpha_2 \nu|_{A_2} + \alpha_3 \nu|_{A_3}$$

will correspond to an element of $C$ provided that $\alpha_i \geq 0$ and

$$\alpha_1 a_1 + \alpha_2 a_2 + \alpha_3 a_3 = 1,$$
$$\alpha_1 b_1 + \alpha_2 b_2 + \alpha_3 b_3 = v.$$

Clearly, this system of linear equations is solved by the vector $(1,1,1)$ but admits also another, different solution $(\alpha_1, \alpha_2, \alpha_3)$ with $0 \leq \alpha_i \leq 2$. But then $\gamma_i := 2 - \alpha_i$ defines yet another solution. Letting

$$\nu_2 := \gamma_1 \nu|_{A_1} + \gamma_2 \nu|_{A_2} + \gamma_3 \nu|_{A_3},$$

we have found two measures $\nu_1, \nu_2$ corresponding to two functions $f_1, f_2$ in $C$ such that

$$\nu = \tfrac{1}{2}(\nu_1 + \nu_2)$$

and, hence,

$$f = \tfrac{1}{2}(f_1 + f_2).$$

Thus, $f$ cannot be an extreme point of $C$.

Next, we turn to the case in which $f \in C$ satisfies $\sup_x f(x) < 1$. Then the corresponding measure $\nu$ is a true subprobability measure: $\nu(S) < 1$. If $f$ takes more than one value in $(0,1)$, then $\nu$ puts positive charge on two disjoint sets $A_1, A_2$, of which we may assume that $A_1 \cup A_2 = S$. Letting again $a_i := \nu(A_i) > 0$ and $b_i := \int_{A_i} G\,d\nu$, we see that $\nu_1 := \alpha_1 \nu|_{A_1} + \alpha_2 \nu|_{A_2}$ will correspond to some function $f_1 \in C$ provided that $\alpha_i \geq 0$ and

$$\alpha_1 a_1 + \alpha_2 a_2 = m,$$
(39)
$$\alpha_1 b_1 + \alpha_2 b_2 = v,$$

where $m$ may be any number between 0 and 1. The vector $(1,1)$ solves (39) for $m = \nu(S) < 1$. If (39) admits also other nonnegative solutions, then we can argue as in the case $\nu(S) = 1$ that $f$ is not an extreme point of $C$. If the solution to (39) is unique, we take $\varepsilon > 0$ such that the numbers $m_\pm := \nu(S) \pm \varepsilon$ belong to $[0,1]$ and such that the solutions $(\alpha_1^\pm, \alpha_2^\pm)$ corresponding to $m_\pm$ have nonnegative components. Then the measures $\nu_\pm := \alpha_1^\pm \nu|_{A_1} + \alpha_2^\pm \nu|_{A_2}$ correspond to functions $f_\pm \in C$ such that $f = \frac{1}{2}(f_+ + f_-)$, and so $f$ is not an extreme point of $C$.

Finally, consider a function $f \in C$ of the form $f = \beta \mathbb{I}_{J_0} + \mathbb{I}_{J_1}$, where $\beta \in (0,1)$ but the interval $J_0$ is a $\mu$-nullset. In this case, $\beta$ can be changed arbitrarily without violating the condition $\int f \, d\mu = v$, and so $f$ cannot be an extreme point of $C$. $\square$

LEMMA 5.2. *The set $C$ admits an integral representation with respect to its extreme points*: *For every function $f_0 \in C$, there exists a probability measure $\eta$ on $C \cap \mathcal{J}$, defined on the $\sigma$-algebra generated by the maps $f \mapsto f(x)$, $x \in (0, \infty)$, such that*

$$f_0 = \int_{C \cap \mathcal{J}} f \eta(df).$$

PROOF. Consider the affine coding of a function $f \in C$ by a subprobability measure $\nu$ on $S$ as introduced in the proof of Lemma 5.1. By adding an additional point $\partial$ to $S$, we can uniquely extend $\nu$ to a probability measure on $\overline{S} := S \cup \{\partial\}$. The function $G$ defined in (38) will be extended to $\overline{S}$ by letting $G(\partial) := 0$. Then $C$ can be identified with the set $H$ of all Borel probability measures $\nu$ on $\overline{S}$ such that $\int G \, d\nu = v$. Corollary 3 of [24] states that $H$ enjoys an integral representation, which then carries over to $C$ by means of Fubini's theorem. $\square$

PROOF OF PARTS (b) AND (d) OF THEOREM 4.1. In proving (b), our task is to minimize $E[g_k(\varphi)f(\varphi)] = \int gf \, d\mu$ over the set $C$, where $g(x) = g_k(x)/x$. Let $f_0$ be a minimizer in $C$ [which must exist, e.g., by Proposition 2.5 and part (c) of Theorem 4.1] and consider the integral representation $f_0 = \int_{C \cap \mathcal{J}} f \eta(df)$ of Lemma 5.2. Then, according to Fubini's theorem,

$$\int gf_0 \, d\mu \geq \inf_{f \in \text{supp}\,\eta} \int gf \, d\mu,$$

so that there must also be a minimizer in $C \cap \mathcal{J}$.

(d) According to the argument in the proof of part (b), uniqueness of solutions in $C \cap \mathcal{J}$ implies uniqueness in $C$. Moreover, by Proposition 2.5, every $\sigma(\varphi)$-measurable solution $X^*$ is of the form $X^* = f^*(\varphi)$ for some function $f^* \in C$. $\square$

**Acknowledgment.** The author expresses his thanks to the Mathematics Department of the University of British Columbia for the hospitality during the time when this research was carried out.

## REFERENCES

[1] ARTZNER, P., DELBAEN, F., EBER, J.-M. and HEATH, D. (1999). Coherent measures of risk. *Math. Finance* **9** 203–228. MR1850791
[2] BEDNARSKI, T. (1981). On solutions of minimax test problems for special capacities. *Z. Wahrsch. Verw. Gebiete* **58** 397–405. MR639148
[3] CVITANIĆ, J. (2000). Minimizing expected loss of hedging in incomplete and constrained markets. *SIAM J. Control Optim.* **38** 1050–1066.
[4] CVITANIĆ, J. and KARATZAS, I. (2001). Generalized Neyman–Pearson lemma via convex duality. *Bernoulli* **7** 79–97. MR1811745
[5] DANA, R.-A. and CARLIER, G. (2002). Core of convex distortions of a probability on an non atomic space. Preprint, Ceremade, Univ. Paris-Dauphine.
[6] DELBAEN, F. (2002). Coherent measures of risk on general probability spaces. In *Advances in Finance and Stochastics. Essays in Honour of Dieter Sondermann* (K. Sandmann and P. J. Schönbucher, eds.) 1–37. Springer, Berlin. MR1929367
[7] DELBAEN, F. (2000). *Coherent Risk Measures.* Scuola Normale Superiore, Classe di Scienze, Pisa.
[8] DENNEBERG, D. (1994). *Non-additive Measure and Integral.* Kluwer Academic, Dordrecht.
[9] FÖLLMER, H. and LEUKERT, P. (1999). Quantile hedging. *Finance Stoch.* **3** 251–273. MR1842286
[10] FÖLLMER, H. and LEUKERT, P. (2000). Efficient hedging: Cost versus shortfall risk. *Finance Stoch.* **4** 117–146. MR1780323
[11] FÖLLMER, H. and SCHIED, A. (2002). Convex measures of risk and trading constraints. *Finance Stoch.* **6**.
[12] FÖLLMER, H. and SCHIED, A. (2002). Robust representation of convex measures of risk. In *Advances in Finance and Stochastics. Essays in Honour of Dieter Sondermann* (K. Sandmann and P. J. Schönbucher, eds.) 39–56. Springer, Berlin.
[13] FÖLLMER, H. and SCHIED, A. (2002). *Stochastic Finance*: *An Introduction in Discrete Time.* Springer, Berlin.
[14] GILBOA, I. and SCHMEIDLER, D. (1989). Maxmin expected utility with non-unique prior. *J. Math. Econ.* **18** 141–153.
[15] HUBER, P. (1981). *Robust Statistics.* Wiley, New York.
[16] HUBER, P. and STRASSEN, V. (1973). Minimax tests and the Neyman–Pearson lemma for capacities. *Ann. Statist.* **1** 251–263. MR356306
[17] KIRCH, M. (2002). Maximin-optimal tests and least favorable pairs for concave power functions. Preprint, TU Wien.
[18] KIRCH, M. (2002). Efficient hedging in incomplete markets under model uncertainty. Preprint, TU Wien.

[19] KULLDORFF, M. (1993). Optimal control of favorable games with a time limit. *SIAM J. Control Optim.* **31** 52–69.
[20] KUSUOKA, S. (2001). On law invariant coherent risk measures. *Adv. Math. Econ.* **3** 83–95.
[21] ÖSTERREICHER, F. (1978). On the construction of least favourable pairs of distributions. *Z. Wahrsch. Verw. Gebiete* **43** 49–55.
[22] PHAM, H. (2002). Minimizing shortfall risk and applications to finance and insurance problems. *Ann. Appl. Probab.* **12** 143–172.
[23] RIEDER, H. (1977). Least favourable pairs for special capacities. *Ann. Statist.* **5** 909–921.
[24] VON WEIZSÄCKER, H. and WINKLER, G. (1979). Integral representations in the set of solutions of a generalized moment problem. *Math. Ann.* **246** 23–32. MR554129
[25] WINKLER, G. (1988). Extreme points of moment sets. *Math. Oper. Res.* **13** 581–587.

INSTITUT FÜR MATHEMATIK, MA 7-4
TECHNISCHE UNIVERSITÄT BERLIN
STRASSE DES 17. JUNI 136
10623 BERLIN
GERMANY
E-MAIL: schied@math.tu-berlin.de